\documentclass[12pt,a4paper]{amsart}
\usepackage{amsmath}\usepackage{amssymb}
\usepackage{amsthm}
\usepackage[english]{babel}\usepackage{german}
\usepackage{nicefrac}  \usepackage{floatflt}
\usepackage{graphicx}    
\usepackage{xcolor}
\graphicspath{ {./Images/} }
\makeatletter
\newtheorem*{rep@theorem}{\rep@title}
\newcommand{\newreptheorem}[2]{%
\newenvironment{rep#1}[1]{%
 \def\rep@title{#2 \ref{##1}}%
 \begin{rep@theorem}}%
 {\end{rep@theorem}}}
\makeatother

\usepackage{biblatex}
\author{Lauritz Streck}
\date{\today}
\title{Non-Concentration of Primes in $\Gamma \backslash PSL_2(\mathbb{R})$}

\begin{document}
\newtheorem{theorem}{Theorem}[section]
\newreptheorem{theorem}{Theorem}
\newtheorem*{theorem*}{Theorem}
\newtheorem{definition}[theorem]{Definition}
\newtheorem{remark}[theorem]{Remark}
\newtheorem{lemma}[theorem]{Lemma}
\newreptheorem{lemma}{Lemma}
\newtheorem{corollary}[theorem]{Corollary}
\newtheorem{prop}[theorem]{Proposition}
\newtheorem{claim}[theorem]{Claim}
\newtheorem{observation}[theorem]{Observation}
\newcommand{\eps}{\varepsilon}

\begin{abstract}
This paper generalizes the result of Sarnak and Ubis \cite{sarnak-ubis} about non-concentration of primes in horocycle orbits on $PSL_2(\mathbb{Z}) \backslash PSL_2(\mathbb{R})$ to any lattice in $PSL_2(\mathbb{R})$. The proof combines the asymptotic result of Strömbergsson \parencite{strombergsson} and Venkatesh's method \parencite{venkatesh} with the approach of Sarnak and Ubis of approximating horocycle pieces with periodic horocycles. The key step is to establish a dichotomy between $\{\xi h(t), t \in [0, T] \}$ having good equidistribution in $\Gamma \backslash PSL_2(\mathbb{R})$ and it being approximable by closed horocycle pieces with small period.  In a follow-up paper,  a similar approach will be used to show equidistribution of $\xi h(n^{1+\gamma})$ for small $\gamma>0$, generalizing Venkatesh's result \parencite{venkatesh} to non-compact $\Gamma$.
\end{abstract}
\maketitle
\section{Introduction}
\subsubsection*{General Introduction}
Let $G=PSL_2(\mathbb{R})$ and $\mu_G$ be the Haar measure on $G$. Let $\Gamma$ be a lattice in $G$, that is, a discrete subgroup such that $\mu_X$, the projection of the Haar measure to $X=\Gamma \backslash G$, is finite (and assumed to fulfill $\mu_X(X)=1$). The dynamics of the space $X$ with respect to $\mu_X$ have been studied extensively in recent years, in part because of the strong connection to Diophantine approximation in the case of $\Gamma=PSL_2(\mathbb{Z})$.\\
The group $G$ can be parametrized in terms of
\[
h(x):=\begin{pmatrix}
1 & x \\ 0 & 1
\end{pmatrix} \quad
a(y):=\begin{pmatrix}
y^{\frac{1}{2}} & 0 \\ 0 & y^{-\frac{1}{2}}
\end{pmatrix} \quad k(\theta):=\begin{pmatrix}
\cos \theta & \sin \theta \\ -\sin \theta & \cos \theta
\end{pmatrix},
\]
which induces a natural left-invariant metric $d_G$ on $G$. This metric descends to $X$ via $d_X(\Gamma g, \Gamma h)=\inf_{\gamma \in \Gamma} d_G(g, \gamma h)$.\\
The geodesic flow
\[
g_t(g):=ga (e^t)=\begin{pmatrix}
ae^{\frac{t}{2}} & be^{-\frac{t}{2}} \\ ce^{\frac{t}{2}} & de^{-\frac{t}{2}}
\end{pmatrix}
\]
and the horocycle flow 
\[
h_t(g):=gh(t)=\begin{pmatrix}
a & b+at \\ c & d+ct
\end{pmatrix}
\]
act ergodically on $(X, \mu_X)$. The horocycle orbits were found to exhibit a very rigid behaviour. Fürstenberg showed that $\mu_X$ is uniquely ergodic under $h_t$ when $X$ is compact \parencite{furstenberg}. For a general $\Gamma$, there are periodic orbits under $h_t$, but these carry all other invariant measures. Precisely, Dani and Smillie showed that both $h_t(\xi), t \in \mathbb{R}$ and $h_n(\xi), n \in \mathbb{N}$ equidistribute with respect to $\mu_X$ unless $\xi$ is periodic, in the sense $t \mapsto h_t(\xi)$ periodic \parencite{dani}. As all periodic orbits are isomorphic to the torus, questions are reduced to questions on tori if the point $\xi$ is periodic. \\
With these questions settled, questions about the equidistribution of other orbits were raised. Shah conjectured that $\xi h(n^\alpha), n \in \mathbb{N}$ equidistributes with respect to $\mu_X$ for all $\alpha \geq 1$ and $\xi$ non-periodic \cite{shah}. Margulis conjectured that $\xi h(p)$ would equidistribute with respect to $\mu_X$, where $p$ is running over the primes and $\xi$ is non-periodic \cite{margulis}. This paper provides partial progress in the latter question by proving that primes do not concentrate anywhere.\\
The way to showing non-concentration is through controlling averages of the form $\frac{s}{T}\sum_{sn \leq T} f(ph(sn))$ and applying sieve methods.  One natural way to do this is to use a smooth approximation of the primes and show equidistribution of this object,  which we will do in this paper.  We will take the pseudo random measure $\nu$,  which gets introduced below. Hopefully,  the reader will find this presentation more coherent and easier to generalize.  Using the Selberg sieve instead like Sarnak and Ubis do in \cite{sarnak-ubis} would also be possible.  \\
Green and Tao used their pseudo random measure $\nu$ as a smooth approximation of the primes to prove their celebrated theorem of the prime numbers containing arbitrarily long arithmetic progressions \cite{green-tao}. Like in the case of proving properties in function spaces through approximation by smooth functions, introducing $\nu$ allowed them to split the proof of some properties of the primes into two parts: First, showing that $\nu$ has certain properties (like the pseudo randomness condition of the $k$-linear forms condition in \cite{green-tao}) and second, that one can recover properties of the primes through these properties of $\nu$ (the relative Szemer\'edi theorem in the case of arithmetic progressions). We will show in this paper that horocycle orbits along~$\nu$ equidistribute. \\
Goldston and Yildirim defined 
\[
\Lambda_R(n):=\sum_{k < R, \; k|n} \mu(k) \log\left(\frac{R}{k}\right),
\]
modeled after the standard convolution identity for the von Mangoldt-function \cite{goldston}. Green and Tao \cite{green-tao} defined the pseudo random measure $\nu$ by
\[
\nu(n):=\frac{1}{\log R} \Lambda_R^2(n).
\]
Due to the particularities of finding arithmetic sequences,  they restricted $\nu$ to a coprime residue class of some big integer $W$,  which is commonly called the $W$-trick.  In our setting, it will not be necessary.\\
Define furthermore $\mathrm{dist}(\Gamma g):=d_X(\Gamma g, p_0)$ where $p_0 \in X$ could be any point, for example $p_0=\Gamma$ (for the definition of $d_X$, see Section 3). For a function $f \in C^4(X)$ let $\Vert f \Vert_{W^4}$ be its Sobolev norm in the Hilbert space $W^{4, 2}$ involving the fourth derivative and let $\Vert f \Vert_{\infty, j}$ be the supremum norm of the $j$-th derivatives. Define
\[
\Vert f \Vert:=\Vert f \Vert_{W^4}+\Vert f \Vert_{\infty, 1}+\Vert f \Vert_{\infty, 0}.
\]
The main result of this paper is:
\begin{theorem} \label{main}
Let $\Gamma \subset PSL_2(\mathbb{R})$ be a lattice. Let $R=T^{\theta}$, where $0<\theta \leq  \frac{\beta}{40}$ is fixed.  Here, $\beta$ is the constant depending only on $\Gamma$ appearing in Theorem \ref{venkatesh}.\\
For a non-periodic $\xi \in X$ and a function $f \in C^4(X)$ with $\Vert f \Vert=1$,
\[
\left| \frac{1}{T} \sum_{n \leq T} f(\xi h(n)) \nu(n)-\int f \; d\mu_X \right| \ll r^{-\theta}+\frac{\log \log R}{\log R},
\]
where $r=T \exp(-\mathrm{dist}(g_{\log T}(\xi)))$ and the implied constant depends only on $\Gamma$. Because $r \to \infty$ as $T \to \infty$, the sequence $\xi h(n)$ equidistributes along $\nu$.
\end{theorem}
The biggest possible size for $\theta$
The second error term is due to the normalization of $\nu$, while $r$ rules the equidistribution on the horocycle segment up to time $T$.\\
Theorem \ref{main} recovers some of the properties of the primes. An immediate corollary is the generalization of the result proved by Sarnak and Ubis in \parencite{sarnak-ubis} for $PSL_2(\mathbb{Z})$.
\begin{corollary} \label{nonconcentration}
Take any non-periodic point $\xi \in X$. For any non-negative function $f \in C^4(X)$ with $\Vert f \Vert=1$,
\[
\frac{1}{\pi(T)} \sum_{p \leq T} f(\xi h(p)) \leq \frac{1}{\theta }\int f \; d\mu_X +O\left(r^{-\theta}+\frac{\log \log R}{\log R}\right)
\]
where the sum is over primes, and $r$ and $\theta$ are as in Theorem \ref{main}.\\
In particular, any limit measure of the primes is absolutely continuous with respect to $\mu_X$ and the primes are dense in a set of positive measure.
\end{corollary}
Corollary \ref{nonconcentration} could be proved without using Theorem \ref{main} by adapting the proof of Theorem \ref{main}. One would use sieve methods instead of the normalization of $\nu$ and the Siegel-Walfisz theorem instead of the Siegel-Walfisz type Theorem~\ref{siegelwalfisz} for $\nu$.\\
Unfortunately, these are all properties of the orbit under primes obtainable through Theorem \ref{main}. It falls short both of showing density and equidistribution of primes. The reason is that to recover properties of the primes,  one needs stronger properties of $\nu$ (like the $k$-linear forms condition in \cite{green-tao}, which is stronger than containing arithmetic progressions of length $k$). In the case of equidistribution of primes, one would probably need good equidistribution of $\nu$ not along simple sums, but along sums of the type
\[
\sum_{n \leq \frac{N}{s_1s_2}} f(\xi h(ns_1)) f(\xi h(ns_2)) \nu(n)
\]
for natural numbers $s_1, s_2$. As sums of these types are not well understood at all, even when replacing $\nu$ by $1$, showing equidistriution of primes using $\nu$ seems to require significant new input.\\
The methods employed may be of interest beyond the question in this paper; in particular,  the following result may have other applications.  It will be instrumental in the proof of equidistribution of $\xi h(n^{1+\gamma})$ for small $\gamma$ performed in a follow-up paper.
\begin{lemma} \label{orbitapprox}
Let $p \in X$ and $T \geq 0$. Let $\delta>0$ and $K \leq T$.  There is an interval $I_0 \subset [0,T]$ of size $|I_0| \leq \delta^{-1} K^2$ such that:\\
For all $t_0 \in [0,T] \backslash I_0$,  there is a segment $\{\xi h(t),  t \leq K\}$ of a closed horocycle approximating $\{ph(t_0+t), 0 \leq t \leq K\}$ of order $\delta$,  in the sense that
\[
\forall 0 \leq t \leq K: \quad d_X\left(ph(t_0+t),  \xi h(t)\right) \leq \delta. 
\]
The period $P=P(t_0, p)$ of this closed horocycle is at most $ P \ll r$, where $r=T \exp(-\mathrm{dist}(g_{\log T}(p)))$ is as in Theorem \ref{main}.\\
Moreover,  one can assure $P \gg \eta^2 r$ for some $\eta>0$ by weakening the bound on $I_0$ to $|I_0| \leq \max\left(\delta^{-1} K^2,  \eta T\right)$.
\end{lemma}
This result is useful because it bridges the gap between the compact case with good asymptotics and the periodic case in some sense.

\subsubsection*{Relation to Previous Papers}
Venkatesh showed that for cocompact lattices, $\xi h(ns), 1 \leq n \leq T$,  is equidistributed with error $T^{-\epsilon}$ as long as $s \ll T^\epsilon$ \parencite{venkatesh}. He deduced that $\xi h(n^{1+\gamma})$ equidistributes for sufficiently small $\gamma>0$ (where $\gamma$ and $\epsilon$ depend only on the lattice).  This is the result that will be generalized to all lattices $\Gamma \subset PSL_2(\mathbb{R})$ in the aforementioned follow-up paper.\\
Venkatesh's proof combined the quantitative equidistribution result 
\[
\left| \frac{1}{T} \int_0^T f(\xi h(t)) \; dt - \int f \; d\mu_X \right| \ll_f T^{-2\epsilon} 
\]
for $f \in C^\infty(X)$ (see Lemma 9.4 in \cite{venkatesh}; the ideas go back to Ratner) with a trick to bound the Fourier coefficients. With an argument as in the proof of Proposition \ref{case1}, the theorem of Venkatesh immediately implies Theorem \ref{main} for cocompact $\Gamma$.\\
Strömbergsson proved an effective equidistribution theorem for the horocycle flow in the non-compact case \parencite{strombergsson}. Strömbergsson showed that 
\begin{equation} \label{strombergsson}
\left| \frac{1}{T} \int_0^T f(\xi h(t)) \; dt - \int f \; d\mu_X \right| \ll_f r^{-\alpha},
\end{equation}
where $\alpha$ only depends on $\Gamma$ and $r$ is as in Theorem \ref{main}. \\
This relates the asymptotics of the horocycle flow at time $T$ to the location of $g_{\log T}(\xi)$; the further $g_{\log T}(\xi)$ is up in some cusp, the worse the asymptotics. Strömbergsson's result can be combined with the method of Venkatesh to prove the theorem below, as done for example by Zheng (\cite{zheng}, Theorem 1.2).
\begin{theorem} \label{venkatesh}
Let $\Gamma$ be a non-compact lattice in $G$. Let $f \in C^4(X)$ with $\Vert f \Vert < \infty$ and $1 \leq s<T$. Then
\[
\left| \frac{s}{T} \sum_{1 \leq j \leq \nicefrac{T}{s}} f(\xi h(sj))-\int f \; d\mu_X \right| \ll s^{\frac{1}{2}} r^{-\frac{\beta}{2}} \Vert f \Vert
\]
for any initial point $\xi \in X$, where $r=T \exp(-\mathrm{dist}(g_{\log T}(\xi)))$. The parameter $\frac{1}{6}>\beta>0$ and the implied constant depend only on $\Gamma$.
\end{theorem}
If $r \gg T^\epsilon$ for some $\epsilon>0$, the situation is very similar to the compact case. The set of points $\xi$ with $r \gg T^\epsilon$ for all $T$ has full measure (by geodesic excursion rates, compare the introduction of \cite{strombergsson}). Thus, if one restricts the analysis to a subset of initial points of full measure, as done by Zheng in \parencite{zheng}, Theorem \ref{venkatesh} is all that is needed to show equidistribution of $n^{1+\gamma}$. Similarly, McAdam shows density of almost primes in $SL_n(\mathbb{Z}) \backslash SL_n(\mathbb{R})$ under a Diophantine condition which is equivalent to $r \gg T^\epsilon$ in two dimensions \parencite{mcadam}; compare Remark \ref{diophcond}.\\
Statements about density are not hard to get from Theorem \ref{venkatesh}, because any non-periodic $\xi$ has a sequence $T_i \to \infty$ such that $r(T_i) \gg T_i$. This holds because $g_{\log T}(\xi)$ returns to a compact set infinitely often (compare Lemma \ref{fundperiod}). Explicitly, one immediately gets density of $\xi h(n^{1+\gamma})$ in $X$, density of almost primes in $X$ and density of primes in a set of positive measure (shown as in the proof of Proposition \ref{case1}) from Theorem \ref{venkatesh}.\\
Sarnak and Ubis chose a different approach and analyzed the quantitative equidistribution of $\xi h(sn)$, for \textit{all} $\xi$ and \textit{all} times $T$ \parencite{sarnak-ubis}. They did this in the case $\Gamma=PSL_2(\mathbb{Z})$, defining a fundamental period $y_T$. This fundamental period is based on the imaginary part of the horocycle segment $\xi h([0 ,T])$ and turns out to be closely related to $r$. They then proceeded to show that the horocycle segment $\xi h(t), 0 \leq t \leq T$ is approximable by periodic horocycle with period at most $y_T$.  Analyzing the situation on the periodic horocycle separately, they deduced that for non-negative $f \in C^4$
\[
\frac{1}{\pi(T)} \sum_{p \leq T} f(\xi h(p)) \leq 10 \int f \; d\mu_X+o_T(1)
\]
for all $T$, which implies non-concentration of primes. They did not use Strömbergsson's result and Theorem \ref{venkatesh}, but used estimates of automorphic forms to obtain similar asymptotics for $r \gg T^\epsilon$.
\subsubsection*{Strategy}
We will combine Theorem \ref{venkatesh} with the approach of Sarnak and Ubis to generalize their result to all lattices $\Gamma \subset G$.  The main step is to generalize their fundamental period from $\Gamma=PSL_2(\mathbb{Z})$ to arbitrary $\Gamma$.  This approach culminates in Lemma \ref{orbitapprox}.  This theorem will allow us to reduce the analysis to closed horocycles in the cases when $r \ll T^{\frac{1}{20}}$ and the asymptotics are bad.  On those,  we will use the Siegel-Walfisz type Theorem \ref{siegelwalfisz} for $\nu$ to finish the proof. 

\subsubsection*{Structure of this paper}
Chapter 2 contains the proof of the Siegel-Walfisz theorem for $\nu$ and ends with a short proof of Corollary \ref{nonconcentration} assuming Theorem \ref{main}. Chapter 3 recalls basics of the dynamics on quotients of $G$ and their relation to quotients of the hyperbolic plane.\\
In Chapter \ref{fundamental},  Lemma \ref{orbitapprox} is proven by generalizing the fundamental period $y_T$ to all lattices $\Gamma$,  establishing $r \sim y_T^{-1}$
and analyzing horocycle segments in $PSL_2(\mathbb{R})$. \\
In Chapter 5, Theorem \ref{venkatesh},  Lemma \ref{orbitapprox} and Theorem \ref{siegelwalfisz} are combined to prove Theorem \ref{main}.

\subsubsection*{Notation}
Elements and flows in $G=PSL_2(\mathbb{R})$ are denoted by
\[
h(x)=\begin{pmatrix}
1 & x \\ 0 & 1
\end{pmatrix}, \quad
a(y)=\begin{pmatrix}
y^{\frac{1}{2}} & 0 \\ 0 & y^{-\frac{1}{2}}
\end{pmatrix}, \quad k(\theta)=\begin{pmatrix}
\cos \theta & \sin \theta \\ -\sin \theta & \cos \theta
\end{pmatrix},
\]
$g_t(g)=ga(e^t)$ and $h_t(g)=g h(t)$. We set $e(t):=\exp(2\pi i t).$\\
A fundamental domain of $X=\Gamma \backslash G$ is chosen to be the unit tangent bundle $E=T_1 F$ of a connected fundamental domain $F$ of the action of $\Gamma$ on the upper half plane $\mathbb{H}$. The open interior is denoted by $F^\mathrm{o}$ and the closure by $\overline{F}$. For $g \in G$, $g.z=\frac{az+b}{cz+d}$, and $g.i$ is the projection to $\mathbb{H}$. When we write $\mathrm{Im}(g)$, we mean $\mathrm{Im}(g.i)$.\\
The objects $r$, $\mathrm{dist}(\xi)$ and the norm $\Vert f \Vert$ of $f \in C^4(X)$ are defined before Theorem \ref{main}. The definition of $y_T$ can be found in Chapter \ref{fundamental} on Page \pageref{fundperdef}. \\
The inequality $f \ll g$ and $f=O(g)$ mean that there is an absolute constant $C$ such that $|f(t)| \leq C |g(t)|$ for all $t$. Write  $f \sim g$ if $f \ll g$ and $g \ll f$. Unless equipped with an index to indicate further dependence, the implicit constants only depend on $\Gamma$.\\
The divisor  function $\tau(n)$ counts the divisors of $n$ and the Euler totient function is denoted by $\phi(n)$. The Möbius function $\mu$ is supported on square free numbers and is defined by
\[
\mu(1)=1, \quad \mu(p_1 \dots p_n)=(-1)^n
\]
for distinct prime numbers $p_1, \dots, p_n$. Its square $\mu^2$ is the characteristic function of square free numbers. For integers $e$ and $d$, the least common multiple is denoted by $[e, d]$ and the greatest common divisor by$(e, d)$.

\subsubsection*{Acknowledgements}
This text was written as my master's thesis at the  Hebrew University in Jerusalem. I am very grateful for the warm welcome, the superb learning environment and the generous support I received at the Hebrew University, especially from Hillel Fürstenberg, Elon Lindenstrauss, Shahar Mozes, Jane Turner and my advisor Tamar Ziegler.\\
Initially, I came to Jerusalem in an exchange year. After an amazing year, I decided to finish my entire degree in Jerusalem.  Many thanks to the entire department and especially to Elon and Tamar, who organized support me in my second year. I am thankful for the many discussions with Tamar about the project and to Elon's open ear for questions whenever I had any.\\
One of the most astonishing experiences in this year for me was a 1on1 reading course with Hillel about his proof of Szemer\'edi's theorem. I learned a lot from his explanations of the proof and his views on the life of a mathematician in general, patiently delivered in his small office amidst piles of books and theses, veiled in layers of chalk dust.  I am very happy that I came to a university where even the most senior researchers (Hillel was 83 at the time) still come to their offices and gladly teach young students. \\
A big thank you to the anonymous referee, whose review improved this paper tremendously.  The suggested tweaks to the proofs made the ideas much clearer and saved close to 10 pages of calculation.  If the reader finds the proofs intuitive,  the changes made after the review have a substantial part in that. \\
Moreover,  the review did not only help the reader. It also made the ideas clearer to me and got me thinking about the material again. This led me to a proof of equidistribution of $ph(n^{1+\gamma})$,  which will be the subject of a follow-up paper. Without the very helpful review,  this would in all likelihood not have happened.\\
Lastly,  I want to thank Dan Goldston, Andreas Strömbergsson, Peter Sarnak and especially Adrián Ubis for their helpful replies to questions I asked.
\newpage

\section{Properties of $\nu$}
In this section,  we are going to derive the Siegel-Walfisz type Theorem \ref{siegelwalfisz} for $\nu$ and prove that Theorem \ref{main} implies Corollary \ref{nonconcentration}.
\begin{lemma} \label{lemma1}(Lemma 2.1 in \cite{goldston})
Let $R>0, k \in \mathbb{N}$ such that $\log(k) \ll \log(R)$.  Then
\[
\sum_{d \leq R, (d, k)=1} \frac{\mu(d)}{d} \log\left(\frac{R}{d}\right)=\frac{k}{\phi(k)}+O\left(\frac{1}{\log^3(R)}\right)
\]
and
\[
\sum_{d \leq R, (d, k)=1} \frac{\mu(d)}{\phi(d)} \log\left(\frac{R}{d}\right)=\mathfrak{S}_2(k)+O\left(\frac{1}{\log^3(R)}\right),
\]
where $\mathfrak{S}_2$ is the singular series from the Goldbach conjecture, supported on positive even numbers and given by $\mathfrak{S}_2(2n)= 2 C_2 \prod_{p|n, p>2} \left(\frac{p-1}{p-2}\right)$ with $C_2=\prod_{p>2} \left(1-\frac{1}{(p-1)^2} \right)$. 
\end{lemma}
The next lemma from \cite{goldston} is cited here only in the case $j=1$ with simplified error terms. The validity of the simplifications can be seen from their remarks after Lemma 2.2 and the fact that $\mathfrak{S}_2(k) \ll \tau(k)$. 
\begin{lemma} \label{lemma2} (Lemma 2.4 in \cite{goldston})
Let $R\geq 1$ and $k \in \mathbb{N}$ such that $\log k \ll \log R$. Then
\[
\sum_{d \leq R, (d, k)=1} \frac{\mu^2(d)}{\phi(d)} \mathfrak{S}_2(dk)=\log(R)+O\left(\log \log 3k \right).
\]
\end{lemma}
These lemmas can be combined in a similar way to the proof of Theorem 5.1 in \cite{goldston} to yield the following proposition. 
\begin{prop} \label{nuestimate}
Let $I$ be an interval in $\mathbb{N}$ and $R>1$. Let $q \in \mathbb{N}$ such that $\log q \ll \log R$ and let $j \leq q$ be coprime to $q$.  Then
\begin{align*}
\frac{1}{|I|}\sum_{n \in I} \Lambda_R^2(qn+j)&=\frac{q}{\phi(q)} \log R+\frac{q }{\phi(q)} O\left( \log \log R  \right)+O\left(\frac{R^2}{|I|}\right).
\end{align*}
\end{prop} 
\begin{proof}
Note that for any integer $k$, 
\[
\sum_{\substack{n \in I \\ k | qn+j  }} 1=\begin{cases} \frac{|I|}{k}+O(1),  \; \; (q, k)=1  \\ 0, \;\; \mathrm{else} \end{cases}.
\]
To bound the sums of the appearing error terms, we record for later use that
\begin{equation} \label{firsterror}
\sum_{k \leq R} \frac{1}{k \log\left(R/k\right)} \leq \sum_{j \leq \log R} \sum_{\frac{R}{2^{j+1}} \leq k \leq \frac{R}{2^{j}}} \frac{1}{k \log\left(2^j\right)} =O(\log \log R)
\end{equation} 
and, similarly, using the well known bound $\frac{1}{\phi(k)} \ll \frac{\log \log k}{k}$,
\begin{equation} \label{seconderror}
\sum_{k \leq R} \frac{1}{\phi(k) \log^2\left(R/k\right)} =O(\log \log R).
\end{equation}
So let us bound the terms.  Unboxing the definition, we see that
\begin{align*}
\sum_{n \in I} \Lambda_R^2(qn+j)&=\sum_{d, e \leq R} \mu(d) \mu(e) \log\left(\frac{R}{d}\right) \log\left(\frac{R}{e}\right) \sum_{\substack{n \in I \\ d, e | qn+j  }} 1\\
&=O(R^2)+|I|  \sum_{\substack{d, e \leq R\\ (de, q)=1}}\frac{ \mu(d) \mu(e) \log\left(\frac{R}{d}\right) \log\left(\frac{R}{e}\right)}{[d, e]}
\end{align*}
where $[d, e]$ is the least common multiple. Denote by $\sum^{\prime}$ a sum in which all summation variables are coprime to $q$ and to each other. \
We find by using Lemma \ref{lemma1} that  
\begin{align*}
&\sum_{\substack{d, e \leq R\\ (de, q)=1}}\frac{ \mu(d) \mu(e) \log\left(\frac{R}{d}\right) \log\left(\frac{R}{e}\right)}{[d, e]}\\&=\sum^{\prime}_{\substack{m \leq R\\md, me \leq R}} \frac{\mu^2(m)}{m} \frac{\mu(d)}{d} \frac{\mu(e)}{e} \log\left(\frac{R}{md}\right) \log\left(\frac{R}{me}\right)\\
&=\sum^{\prime}_{md \leq R} \frac{\mu^2(m)}{m} \frac{\mu(d)}{d} \log\left(\frac{R}{md}\right) \sum_{\substack{e \leq \nicefrac{R}{m} \\ (e, mdq)=1}} \frac{\mu(e)}{e} \log\left(\frac{\nicefrac{R}{m}}{e}\right)\\
&=\frac{q}{\phi(q)} \left(\sum^{\prime}_{md \leq R} \frac{\mu^2(m)}{\phi(m)} \frac{\mu(d)}{\phi(d)} \log\left(\frac{R}{md}\right) \right)+O(\log \log R),
\end{align*}
where (\ref{firsterror}) was used to bound the error sum coming from Lemma \ref{lemma1}.\\
Applying Lemma \ref{lemma1}, to the main term, we find
\begin{align*}
&\sum^{\prime}_{md \leq R} \frac{\mu^2(m)}{\phi(m)} \frac{\mu(d)}{\phi(d)} \log\left(\frac{R}{md}\right)=\sum_{\substack{m \leq R \\ (m, q)=1}} \frac{\mu^2(m)}{\phi(m)} \sum_{\substack{d \leq \nicefrac{R}{m} \\ (d, mq)=1}} \frac{\mu(d)}{\phi(d)} \log\left(\frac{\nicefrac{R}{m}}{d}\right)\\
&=\sum_{\substack{m \leq R \\ (m, q)=1}} \frac{\mu^2(m)}{\phi(m)} \mathfrak{S}_2(mq)+O\left(\log \log R \right),
\end{align*}
where the error term was bounded with the help of (\ref{seconderror}).  
Applying Lemma \ref{lemma2} finishes the proof.
\end{proof}
\begin{theorem}[Siegel-Walfisz for $\nu$] \label{siegelwalfisz} 
Let $N \in \mathbb{N}$ and let $\nu(n)=\frac{\Lambda_R^2(n)}{\log R}$ with $R$ of size $N^\epsilon$. For any $q \in \mathbb{N}$,  any interval $I \subset[0, N]$ of size $|I| \geq qR^3$ and any $q$-periodic function $f$, 
\[
\frac{1}{|I|} \sum_{n \in I} f(n) \nu(n)=\frac{q}{\phi(q)|I|}\sum_{\substack{n \in I \\ (n, q)=1}} f(n)+O\left(\frac{ \Vert f\Vert_\infty \log \log R}{\log R}\right),
\]
where $\Vert f \Vert_\infty=\max_{r \leq q} |f(r)|$.
In particular,
\[
\frac{1}{N} \sum_{n \leq N} \nu(n)=1+O\left(\frac{ \log \log R}{\log R}\right)
\]
\end{theorem} 
\begin{proof}
Without loss of generality, assume that $\Vert f\Vert_\infty =1$. Fix an interval $I$. For $q=1$, by Proposition \ref{nuestimate},
\[
\frac{1}{|I|} \sum_{n \in I} \nu(n)=1+O\left(\frac{ \log \log R}{\log R}\right)+O\left(\frac{1}{R}\right).
\] 
For general $q \leq N$ and some coprime $j < q$, Proposition \ref{nuestimate} applied to $q^{-1} (I-j)$ (the assumption $\log q \ll \log R$ is satisfied because $R$ is of size $N^{\eps}$) implies that 
\[
\frac{q}{|I|}\sum_{\substack{n\in I\\n\equiv j (q)}} \nu(n)=\frac{q}{\phi(q)}+\frac{q}{\phi(q)} O\left(\frac{\log \log R}{\log R}\right)+O\left(\frac{1}{R}\right).
\]
In particular,  all residue classes coprime to $q$ contribute the same amount to the sum;  it is thus only left to check that the residue classes which are not coprime have negligible contribution. To see this, average over all residue classes $j$ coprime to $q$,  giving
\[
\frac{1}{|I|}\sum_{\substack{n \in I\\ (n, q)=1}} \nu(n)=1+O\left(\frac{\log \log R}{\log R}\right)+O\left(\frac{1}{R}\right),
\]
and compare with the contribution of all residue classes above.  
\end{proof}
To conclude, we will prove Corollary \ref{nonconcentration} assuming Theorem \ref{main}.
\begin{proof}
By partial summation and the prime number theorem, 
\[
\frac{1}{\pi(T)} \sum_{p \leq T} f(\xi h(p))=\frac{1}{T} \sum_{n \leq T} f(\xi h(n)) \;  \tilde{\Lambda}(n)+O\left(\frac{\Vert f \Vert_\infty}{\log T} \right),
\]
where $\tilde{\Lambda}$ is given by
\[
\tilde{\Lambda}(n)=\begin{cases} \log(p), \; n=p  \; \mathrm{prime}\\ 0, \; \mathrm{else} \end{cases}.
\]
By definition of $\nu$, for all primes $p>R$
\[
\nu(p)=\log R=\frac{1}{\theta} \log T.
\]
In particular, $\tilde{\Lambda}(n) \leq \frac{1}{\theta} \nu(n)$ on the interval $[T^\theta, T]$. The corollary follows from Theorem \ref{main}.
\end{proof}

\newpage

\section{Basics of the Dynamics on Quotients of the Hyperbolic Plane}
In this section,  we recall some basics of the hyperbolic plane.  
Any matrix $g \in G$ can be uniquely written as
\[
g=h(x)a(y)k(\theta)
\]
for $x \in \mathbb{R}, y>0, \theta \in [-\nicefrac{\pi}{2}, \nicefrac{\pi}{2})$; this is the Iwasawa parametrization. Moreover, $G$ has a natural left-invariant metric given by
\[
d_G t^2 =\frac{dx^2+dy^2}{y^2}+d\theta^2.
\]
The upper half plane $\mathbb{H}$ carries the hyperbolic metric, which is invariant under the action of $G$ via Möbius transformations $z \mapsto g.z$. This action lifts to an action on the unit tangent bundle $T_1\mathbb{H}$ given by 
\[
\begin{pmatrix}
a & b \\ c & d
\end{pmatrix}.(z, v)=\left(\frac{az+b}{cz+d}, \frac{v}{(cz+d)^2} \right).
\]
There is a natural bijection $PSL_2(\mathbb{R}) \to T_1\mathbb{H}$ given by 
\[
h(x)a(y)k(\theta) \mapsto \left(x+iy, e^{i2\theta}\right)
\]
which is isometric with respect to the metric on $T_1\mathbb{H}$ induced by the hyperbolic metric on $\mathbb{H}$. The measure 
\[
\mu_G=\frac{dx dy}{y^2} \frac{d\theta}{\pi}
\]
is a Haar measure on $G$ ($G$ is unimodular, so there is no distinction between a left and right Haar measure). \\
Fix a lattice $\Gamma$ in $G$ and set $X:=\Gamma \backslash G$. This space carries a left-invariant metric 
\[
d_X(\Gamma g_1, \Gamma g_2)=\min_{\gamma \in \Gamma} d_G(g_1, \gamma g_2)
\]
and a finite measure $\mu:=\mu_X=\pi_\# \mu_G$, the push forward of $\mu_G$ under the projection to $X$. Write $p.i$ for $\Gamma g.i \in \Gamma \backslash \mathbb{H}$ with $p=\Gamma g \in X$.\\
The fundamental domain $E$ of $X$ in $G$ can be chosen to be the tangent bundle of a convex hyperbolic polygon $F \subset \mathbb{H}$ (i. e. a polygon in which each edge is a geodesic) with finitely many edges and vertices; this $F$ is a fundamental domain for $\Gamma \backslash \mathbb{H}$. Explicitly, one can choose $F$ to be a Dirichlet domain for any point $z_0 \in \mathbb{H}$; that is, the interior of $F$ is given by 
\[
F^{\mathrm{o}}=D(z)=\{z \in \mathbb{H}| d_G(z, z_0) < d_G(z, \gamma z_0) \; \forall \gamma \in \Gamma\backslash \{\mathrm{id}\}\}.
\]
The polygon $F$ might or might not have boundary vertices, i. e. point of adjacency with $\partial \mathbb{H}=\mathbb{R} \cup \infty$. $X$ is compact if and only if there are none of those. If there are some, the equivalence class of each boundary vertex with respect to $\Gamma$ is called a cusp of $X$ (for example if $\Gamma=PSL_2(\mathbb{Z})$, $X$ has the single cusp $\infty$ which is equivalent to every rational number). For each representative $r_i \in \partial \mathbb{H}$ of each cusp $\Gamma r_i$ there is a fundamental domain $F$ such that all other boundary vertices of $F$ are inequivalent to $r_i$; this is because we can take some point far up in the cusp as basic point for the Dirichlet domain. Proofs for all of these statements can be found in chapter 11 of \cite{einsiedler-ward}.\\
There are two important flows on $G$. The first one is the geodesic flow given by
\[
g_t(g)=ga ( e^t )=\begin{pmatrix}
ae^{\frac{t}{2}} & be^{-\frac{t}{2}} \\ ce^{\frac{t}{2}} & de^{-\frac{t}{2}}
\end{pmatrix}
\]
and the second one is the horocycle flow given by 
\[
h_t(g)=gh(t)=\begin{pmatrix}
a & b+at \\ c & d+ct
\end{pmatrix}
\]
where $g=\begin{pmatrix}
a & b \\ c & d
\end{pmatrix}$. The flows are well-defined on $X$ because $\Gamma$ acts from the left and the flows act from the right. On $G$ the behaviour of these flows is not very interesting, but on $X$ it is. As outlined in the introduction, the dynamics with respect to the horocycle flow exhibit a very rigid behaviour.\\
There are no periodic horocycle orbits in $X$ if and only if $X$ is compact. If there are periodic orbits, their structure is as follows:
\begin{lemma} \label{perorbits}(Lemma 11.29 in \cite{einsiedler-ward})
Let $\Gamma$ be a lattice such that $X$ is non-compact. To every cusp of $X$ corresponds exactly one one-parameter family of $h$-periodic orbits parametrized by $g_t$; in explanation, if $p \in X$ is such that $t \mapsto ph(t)$ is periodic, then $g_t(p)$ converges to some cusp of $X$ as $t \to \infty$ and all other periodic orbits associated to this cusp contain exactly one element $g_t(p)$ for $t \in \mathbb{R}$. \\
Furthermore, the orbit $t \mapsto ph(t)$ is periodic if and only if $g_t(p) \to \infty$, in the sense that $g_t(p)$ leaves any compact subset of $X$ permanently. 
\end{lemma}
It is shown in the proof (or alternatively, can be seen directly from the statement) that for any boundary vertex of $F$ there is a $\gamma \in \Gamma$ conjugated to $\begin{pmatrix}
1 &1 \\ & 1
\end{pmatrix}$ which fixes this boundary vertex and generates the subgroup of $\Gamma$ fixing the vertex. This $\gamma$ is precisely the one leading to the periodicity of the corresponding orbits. For example for $\Gamma=PSL_2(\mathbb{Z})$ and the sole boundary vertex $\infty$, this matrix is $\gamma=\begin{pmatrix}
1 &1 \\ & 1
\end{pmatrix}$.\\
If $p$ is periodic with period $y$, then $g_t(p)$ is periodic with period $ye^{-t}$ because of the equation $g_t \circ h_s =h_{e^{-t}s} \circ g_t$. 
\newpage

\section{Approximation by closed horocycles} \label{fundamental} 
In this chapter, we will define the fundamental period of a horocycle piece (first defined for $\Gamma=PSL_2(\mathbb{Z})$ in \cite{sarnak-ubis}) and explore the connection to effective equidistribution. The ultimate goal of this section is to prove Lemma \ref{orbitapprox}. \\
Let $n$ be the number of cusps of $X$. Let $r_i \in \partial \mathbb{H}$ be a representative of one cusp of $X$ and let $\gamma_i \in \Gamma$ be the corresponding unipotent element fixing $r_i$ and inducing the periodicity of the corresponding horocycle, as in the discussion after Lemma \ref{perorbits}. Let $\sigma_i \in G$ such that $\sigma_i \gamma_i \sigma_i^{-1}=h(1)$ and $\sigma_i.r_i=\infty$. Explicitly, this $\sigma_i$ consists of a rotation matrix, rotating $r_i$ to $\infty$ and $\gamma_i$ to some $h(t_i)$, and some diagonal element to normalize $t_i=1$.\\
Define $y_i: G \to \mathbb{R}_+$ by 
\[
y_i(g)=\mathrm{Im}(\sigma_i g)
\]
where we mean in slight abuse of notation the imaginary part of the corresponding basepoint $\sigma_i g.i$ in $\mathbb{H}$.  Note that $y_i$ is only well-defined on $G$, not on $X$, and depends on the representative of the cusp.\\
In the natural parametrization $h(x)a(y).i$ for $\mathbb{H}$, every level set with fixed $y$ is a horizontal line, so a circle with point of tangency $\infty$ which is invariant under $h(1)$. We could also parametrize $\mathbb{H}$ by $ \sigma_i^{-1} h(x) a(y) .i$, in terms of which any level set is a circle with point of tangency $r_i$ invariant under $\gamma_i$. $y_i(z)$ then is the $y$-component of $z$ in terms of this new parametrization.

\begin{figure}[h]
\includegraphics[width=\textwidth]{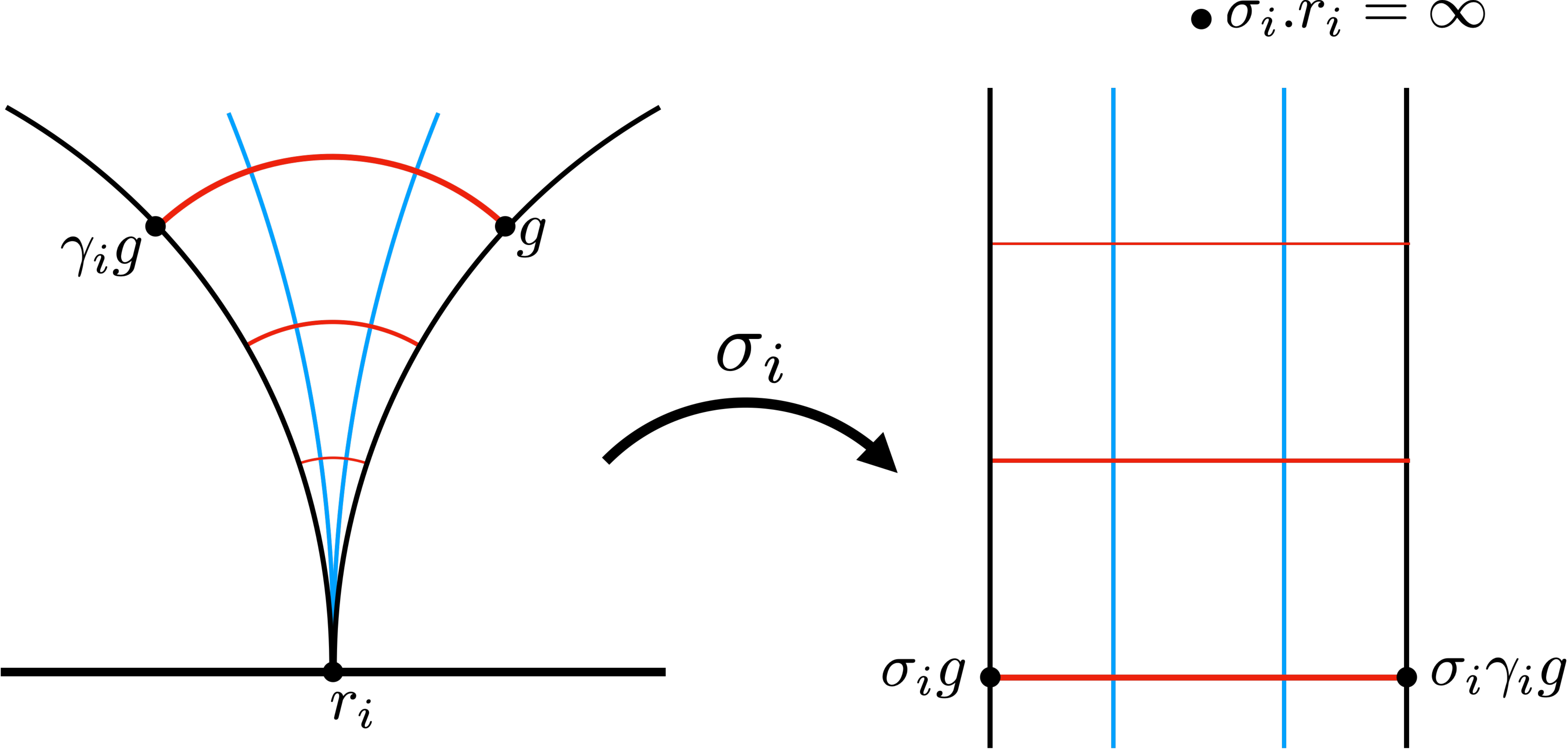}\
\caption{The effects of rotating the Riemann sphere.  If in the new parametrization $\sigma_i g \simeq \textcolor{red}{h(x)} \textcolor{blue}{a(y)}$,  then in the initial parametrization $g \simeq \textcolor{red}{\sigma_i^{-1} h(x) \sigma_i} \textcolor{blue}{\sigma_i^{-1} a(y)}$ and $y_i(g)=y_i(\gamma_i g)=\textcolor{blue}{y}$.  (Here $\simeq$ is defined by $g \simeq h$ if $g$ and $h$ have the same projection to $\mathbb{H}$).}
\end{figure}

Let for $T>0$
\[
Y_i^T(g)=\min\left\{y_i(gh(t))| 0 \leq t \leq T \right\}=\min(y_i(g), y_i(gh(T)))
\]
where the second equality follows from the fact that $\{\sigma_i g h(t)|0 \leq t \leq T\}$ is a piece of a horocycle orbit and thus a segment of a circle in $\mathbb{H}$.\\
Define $y^T_i: X \to \mathbb{R}_+$ by
\[
y^T_i(\Gamma g)=\sup_{\gamma \in \Gamma} Y^T_i(\gamma g).
\]
The supremum is finite and attained for some $\gamma$ because $\Gamma$ is discrete and acts properly discontinuous on $X$. It is independent of the representative of the cusp. Note that $y^0_i: \Gamma \backslash \mathbb{H} \to \mathbb{R}_+$, i. e. $y^0_i$ depends only on the base point $\Gamma g. i$ of $\Gamma g$.\\
Lastly, define
\[  
y_T(p)=\max_{ 1\leq i \leq n}[y_i^T(p)]. \label{fundperdef}
\]
The horocycle piece up to time $T$ starting in $p$ will be close to a periodic horocycle with period $y_T^{-1}$. This is called the \\textit{fundamental period} at time $T$, whose properties we will explore now. 
\begin{lemma}[Parametrization in the Cusps] \label{cusppara}
1. For some small $\epsilon>0$, each $C_i:=\{p| y_i^0(p)^{-1}<\epsilon\} \subset X$ is an open neighborhood of the cusp $\Gamma r_i$ and all $C_i$ are pairwise disjoint. The set $K:=X \backslash \bigcup C_i$ is compact.\\
2. Any point $p$ is contained in $C_i$ if and only if there is a $q$ with the same base point ($p.i=q.i$ in $\Gamma \backslash \mathbb{H}$) such that $qh(t)$ is periodic with period smaller than $\epsilon$. In this case, $qh(t)$ has period $y_i^0(p)^{-1}$.\\
3. Fix disjoint $C_i$ and $K$ as in 1. and take some $p_0 \in K$. Then 
\[
\exp(d_X(p, p_0)) \sim y_0(p)=\max_{1 \leq i \leq n}[y_i^0(p)].
\]
More explicitly, if $p \in C_i$, then $\exp(d_X(p, p_0)) \sim y_i^0(p)$ and if $p \in K$, $\exp(d_X(p, p_0)) \sim 1 \sim y_0(p)$.
\end{lemma}
The implied constants depend on the choice of $\epsilon$, but $C_i$ and $K$ will be fixed from here on.  Parts of the lemma have appeared in the literature before; the function $y_0(p)$ is known as the invariant height function, compare (11) in \cite{strombergsson}.  Part 3. of Lemma \ref{cusppara} was stated as Inequality (14) in the same paper paper,  leaving the proof as an exercise.  We will prove it here for completeness.

\begin{figure}[h]
\includegraphics[width=8cm]{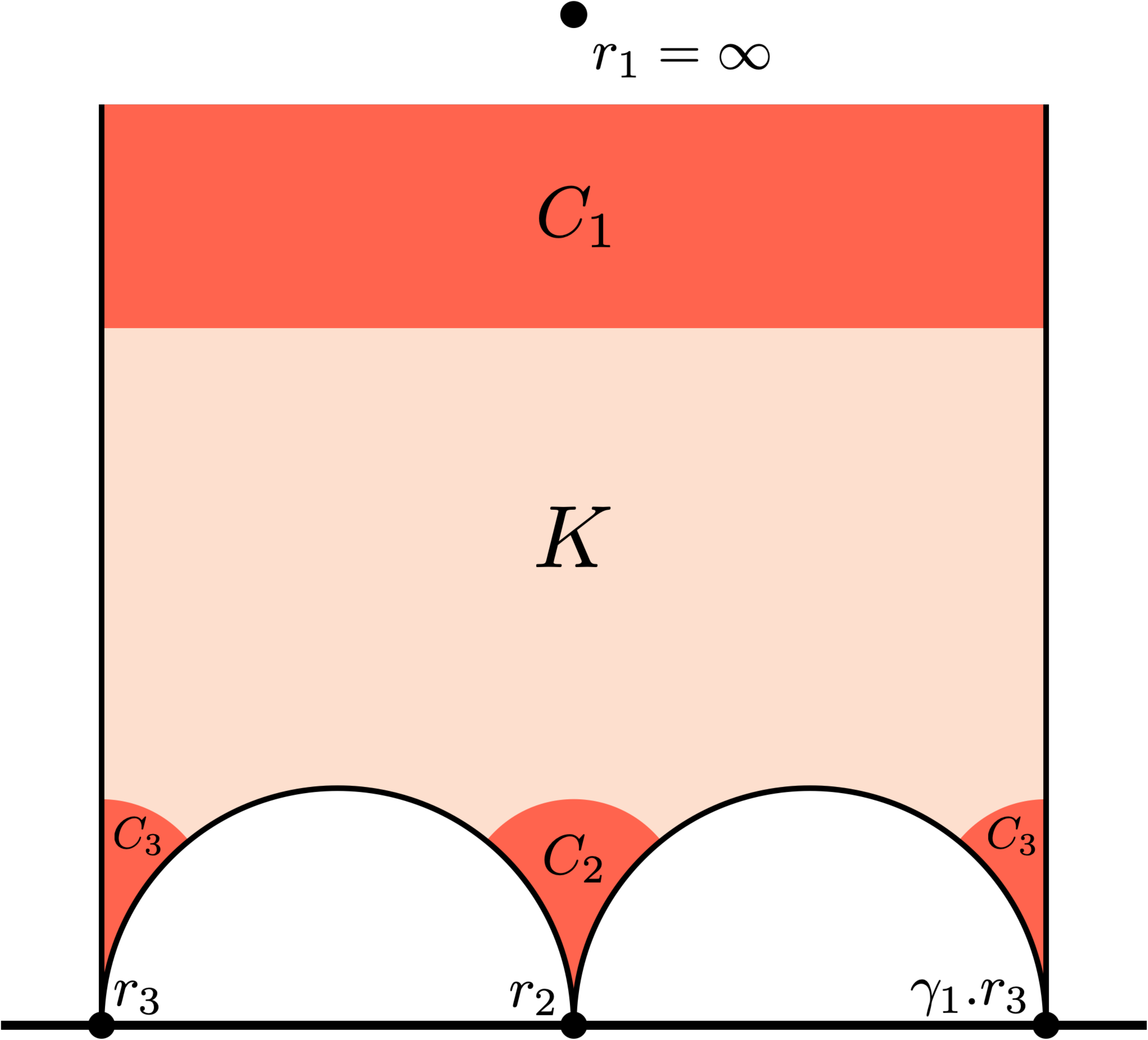}\
\caption{The partitioning of $X$,  drawn on some fundamental domain projected to $\mathbb{H}$.}
\end{figure}

\begin{proof}
As discussed in chapter 3, we can choose a fundamental domain $E$ such that $E=T_1 F$, where the open interior $F^{\mathrm{o}}$ is a Dirichlet domain with boundary vertex $r_i$ such that all other boundary vertices of $F$ are inequivalent to $r_i$. Then $\sigma_i F$ is a fundamental domain of $\sigma_i \Gamma \sigma_i^{-1} \backslash \mathbb{H}$ with boundary vertex $\infty$ which is inequivalent to all other boundary vertices.\\
$\sigma_i F$ is a hyperbolic polygon with finitely many sides. Because $\infty$ is a vertex, two of them must be straight vertical lines. Thus for some big $B \in \mathbb{R}$, $\sigma_i F \cap \{z|\mathrm{Im}(z)>B\}$ is a rectangle. Because there are no other equivalent boundary vertices and $h(1) \in \sigma_i \Gamma \sigma_i^{-1}$ by choice of $\sigma_i$, the horizontal  line of this rectangle has euclidean length $1$ and 
\[
\{z|\mathrm{Im}(z)>B\}=\bigcup_{k \in \mathbb{Z}} h(k) \left(\sigma_i F \cap \{z|\mathrm{Im}(z)>B\} \right).
\]
1. Set $D_i= \Gamma \sigma_i^{-1}  \{z \in \sigma_i F \; |\mathrm{Im}(z)>B\}$ and $C_i=T_1 D_i \subset X$. Let $g \in E$. If $\Gamma g \in C_i$, by definition $y_i^0(\Gamma g) \geq y_i(g)=\mathrm{Im}(\sigma_i g)>B$. If on the other hand $\Gamma g \notin C_i$, because the left translates $h(k)\{z \in \sigma_i F \; |\mathrm{Im}(z)>B\}$ exactly tile the set $\{z|\mathrm{Im}(z)>B\}$, we must have $y_i(\gamma g) \leq B$ for all $\gamma \in \Gamma$ (else $g$ would have two different representatives in the fundamental domain). Thus $y_i^0(g) \leq B$. This shows that $C_i=\{p|y_i^0(p)^{-1} < \epsilon_i\}$ with $\epsilon_i=B^{-1}$.\\
2. Let $p$ now be periodic with period $b < \epsilon$ and let $g$ be the representative of $p$ in $E$. Then $\gamma_i g=gh(b)$ by our choice of $\gamma_i$. The orbit of $\sigma_i  g$ is periodic with respect to infinity in $\sigma \Gamma \sigma^{-1} \backslash G$ because
\[
h(1) \sigma_i g=\sigma_i \gamma_i \sigma_i^{-1} \sigma_i g=\sigma_i g h(b),
\]
But the orbit $U=\{h(x)a(b^{-1})| 0 \leq x \leq 1\}$ is also in $\sigma_i E$ and is periodic with period $b$, so by Lemma \ref{perorbits} they have to agree. Because $b < \epsilon$, we get $\sigma_i p \in U \subset \sigma_i C_i$ and $y_i^0(p)=b^{-1}$. \\
If on the other hand $p=\Gamma g \in C_i$ such that $g \in E$, then $\sigma_i g=(z, v)$ for some $z$ with $\mathrm{Im}(z) > \epsilon^{-1}$. Set $q=\sigma_i^{-1} (z, i)$ and note that $qh(t)$ is periodic with period $\mathrm{Im}(z)^{-1}=y_i^0(p)^{-1}$ because $t \mapsto \sigma_i q h(t)$ is. There cannot be a $q^\prime$ with $q^\prime.i=p.i$ and a smaller period, because $\sigma_i q^\prime. i$ would then have to be at a different level in the fundamental domain $\sigma_i F$.\\
3. If $p \in K$, $\exp(d_X(p, p_0)) \sim 1$ because $K$ is compact. On the other hand, $K.i \subset \Gamma \backslash \mathbb{H}$ is also compact and thus $\overline{K.i} \subset \overline{F}$ as well. Consequently, the continuous functions $y_i$ have to be bounded from above and below on this set, showing $y_0(p) \sim 1$. \\
If $p \in C_i$, $y_i^0(p) > \epsilon^{-1}$ and $y_j^0(p) \leq \epsilon^{-1}$ for $j \neq i$, so $y_0(p)=y_i^0(p)$. We can find a point $\Gamma h \in K$ which is close enough to the cusp $\Gamma r_i$ so that $r_i$ has no equivalent boundary vertices in the Dirichlet domain $D(h.i)$. Consider the fundamental domain $E=T_1 F$ with $F^{\mathrm{o}}=D(h.i)$. Let $g$ be the representative of $p$ in $E$. Note that because $\sigma_i F$ has width at most $1$, $\sigma_i g$ and $ \sigma_i h$ essentially only differ in their $a$ component in the Iwasawa decomposition, that is, their imaginary part. Furthermore, $\mathrm{Im}(\sigma_i h) \sim 1$ because $\Gamma h \in K$.\\
Thus by definition of the Dirichlet domain, the left invariance of the metric and the choice of the fundamental domain,
\begin{align*}
\exp(d_X(p, p_0)) &\sim \exp(d_X(p, \Gamma h))=\exp(d_G(g, h))\\
 &=\exp(d_G(\sigma_i g, \sigma_i h)) \sim \mathrm{Im}(\sigma_i g) =y_i(g)=y_i^0(p).
\end{align*}
\end{proof}
Now we are in the position to establish the connection between the fundamental period and $r$.
\begin{prop} \label{fundperiod} 
Let $p \in X$ and $T \geq 3$. Let $C_i$, $K$ as in Lemma \ref{cusppara} and fix some $p_0 \in K$. Let $r=T e^{-\mathrm{dist}(g_{\log T}(p))}$ be as in Theorem \ref{main}. With $y_T$ defined in the beginning of the chapter,
\[
r^{-1} \sim y_T.
\]
More explicitly, if $g_{\log T}(p) \in C_i$, then $r^{-1}(p) \sim y_i^T(p)$ and if $g_{\log T}(p) \in K$, $r^{-1}(p) \sim T^{-1} \sim y_T(p)$.
All implied constants depend only on the choice of $C_i$ and $p_0$, so ultimately only on $\Gamma$.
\end{prop}

\begin{proof}
In light of Lemma \ref{cusppara}, it suffices to show
\[
y_i^0(g_{\log T}(p)) T^{-1}  \sim y_i^T(p).
\]
Fix a representative $g$ of $p$. Pick some $\gamma \in \Gamma$ and set 
\[
\begin{pmatrix}
a & b \\ c & d
\end{pmatrix}:=\sigma_i \gamma g \in PSL_2(\mathbb{R}),
\]
where we choose $c$ to be non-negative.  Then
\[
Y_i^T(\gamma g)=\min\left(\frac{1}{c^2+d^2}, \frac{1}{c^2+(Tc+d)^2} \right)
\]
by definition; this can be simplified because 
\begin{equation} \label{exercise}
\min\left(\frac{1}{c^2+d^2}, \frac{1}{c^2+(Tc+d)^2} \right) \sim \min\left(\frac{1}{T^2c^2}, \frac{1}{d^2} \right), 
\end{equation}
which is left to the reader as an exercise. The observation that for any $r, s>0$,  
\[
\frac{1}{r+s} \sim \min \left(\frac{1}{r}, \frac{1}{s} \right)
\]
may come in handy showing this. \\
With (\ref{exercise}) in hand,  
\[
T^{-1} Y_i^0(\gamma g_{\log T}(g))=\frac{1}{T^2c^2+d^2} \sim \min\left(\frac{1}{T^2c^2}, \frac{1}{d^2} \right) \sim Y_i^T(\gamma g).
\]
Taking the supremum over $\gamma$ finishes the proof.
\end{proof}
At this point, let us relate two other papers building on the results of Venkatesh and let us translate the respective conditions on the points into our notation.\\
\begin{remark} \label{diophcond}
\upshape In two dimensions the Diophantine condition for $\Gamma g$ (Equation (3.1.c) on page 11 in \cite{mcadam}) of McAdam is 
\[
\min_{\omega \in \mathbb{Z}^2 \backslash \{0\}} \max_{0 \leq t \leq T} \Vert \omega g h(t) \Vert_{\infty}  \gg T^\epsilon 
\]
which translated in the notation of the proof is equivalent to
\[
\min_{\gamma \in SL_2(\mathbb{Z})} \max(|c_\gamma|, |d_\gamma|, |d_\gamma+T c_\gamma|) \gg T^\epsilon
\]
with $\begin{pmatrix}
a_\gamma & b_\gamma \\ c_\gamma & d_\gamma
\end{pmatrix}=\gamma g$. As in Proposition \ref{fundperiod}, the left hand side is asymptotically equal to $y_T^{-\frac{1}{2}}$.\\ 
Zheng proves equidistribution of $n^{1+\eta}$ for points fulfilling a $\kappa$-Diophantine condition, where $\eta$ depends on $\kappa$, \parencite{zheng}. In our notation, this $\kappa=(\kappa_1, \dots, \kappa_n)$-condition for $\kappa_i>0$ says that for any cusp there exist $a_i, b_i >0$ such that either $|c_\gamma|>a_i$ or $|d_\gamma|^{\kappa_i} |c_\gamma|>b_i$ for all $\gamma$. It's easy to check that this condition implies $Y_i^T(\gamma g) \ll T^{-\frac{2}{1+\kappa_i}}$, so by Proposition \ref{fundperiod} again $r \gg T^{\epsilon}$ for $\epsilon<\min_i \frac{2}{1+\kappa_i}$. \\
In both cases Theorem \ref{venkatesh} thus immediately implies the respective results.
\end{remark}
To finish off this section, we prove Lemma \ref{orbitapprox}, giving means to approximate the horocycle segment $ph([0,T])$ with closed horocycles segments of period at most $y_T^{-1} \sim r$.  We restate Lemma \ref{orbitapprox} for the convenience of the reader.
\begin{replemma}{orbitapprox}
Let $p \in X$ and $T \geq 0$. Let $\delta>0$ and $K \leq T$.  There is an interval $I_0 \subset [0,T]$ of size $|I_0| \leq \delta^{-1} K^2$ such that:\\
For all $t_0 \in [0,T] \backslash I_0$,  there is a segment $\{\xi h(t),  t \leq K\}$ of a closed horocycle approximating $\{ph(t_0+t), 0 \leq t \leq K\}$ of order $\delta$,  in the sense that
\[
\forall 0 \leq t \leq K: \quad d_X\left(ph(t_0+t),  \xi h(t)\right) \leq \delta. 
\]
The period $P=P(t_0, p)$ of this closed horocycle is at most $ P \ll r$, where $r=T \exp(-\mathrm{dist}(g_{\log T}(p)))$ is as in Theorem \ref{main}.\\
Moreover,  one can assure $P \gg \eta^2 r$ for some $\eta>0$ by weakening the bound on $I_0$ to $|I_0| \leq \max\left(\delta^{-1} K^2,  \eta T\right)$.
\end{replemma}
\begin{proof}
The quantity $r$ will play no role in the proof; we show that the period of the closed horocycles is bounded by $P \ll y_T^{-1}$ and use Proposition \ref{fundperiod}.  Recall that $y_T$ is a maximum over the different cusps and the elements in $\Gamma$.  Let $\sigma_i$ be the rotation with $\sigma_i.r_i=\infty$ corresponding to the cusp $r_i$ maximizing $y_T$, that is such that $y_T=y_i^T$. The rest of the approximation has nothing to do with $\Gamma$, but is just an observation about approximating horocycle pieces by horizontal lines in $PSL_2(\mathbb{R})$.  The period only comes in because in the coordinate system induced by $\sigma_i$,  the height of the horizontal line is the same as the period of the horocycle in $\Gamma\backslash G$.\\
\begin{figure}[b]
\includegraphics[width=\textwidth]{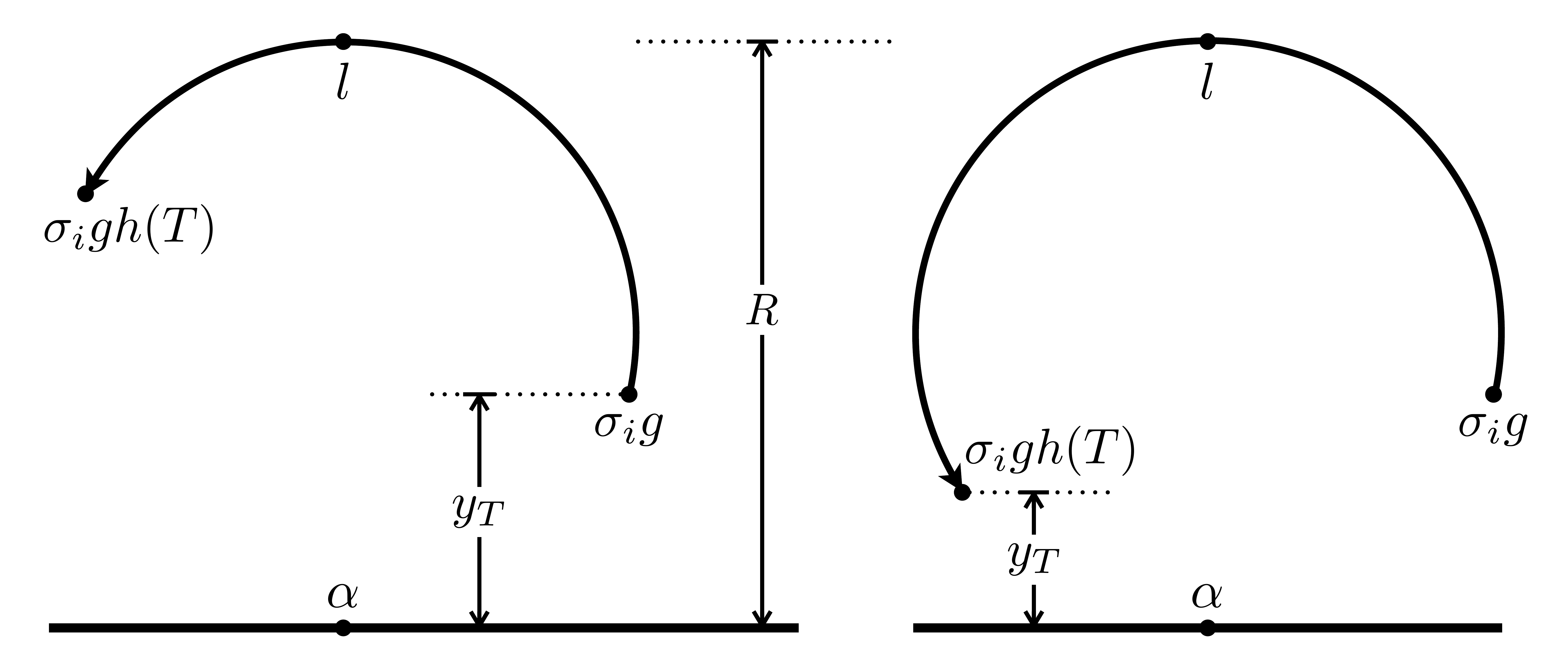} \label{figure3}\
\caption{An overview of the definitions,  to the left when $\sigma_i g$ and to the right when $\sigma_i g h(T)$ minimizes the imaginary part.}
\end{figure}
Let $g$ be a representative of $p$ attaining the supremum in the definition of  $y_i^T$.  The horocycle segment $\{ ph(t), 0\leq t \leq T\}$ is then a circle segment in the modular plane as sketched in Figure \ref{figure3}.  
Write
\[
\sigma_i g=:\begin{pmatrix}
a & b \\ c & d
\end{pmatrix}
\]
and express points on the circle in term of its peak
\[
l:=\sigma_i g h\left(-\frac{d}{c} \right)
=:(\alpha+iR, -i).
\]

In the Iwasawa decomposition,  this is (see (2.3) in \cite{sarnak-ubis})
\[
lh(s)=h\left(\alpha-\frac{Rs}{s^2+1}\right)a\left(\frac{R}{s^2+1}\right)k(-\mathrm{arccot} \; s).
\]
Given some $s$, we will approximate the horocycle segment $\{ lh(s+t), t \leq K\}$ with the periodic horocycle segment  $\{ g_0h(t), t \leq K\}$, where $g_0=:h(x_0)a(y_0)$ lies over the same point in the modular plane as $gh(s)$ and the vector of $g_0$ points straight up.  
The horocycle starting in $g_0$ is then a horizontal line moving right.  

\begin{figure}[t] 
\includegraphics[width=\textwidth]{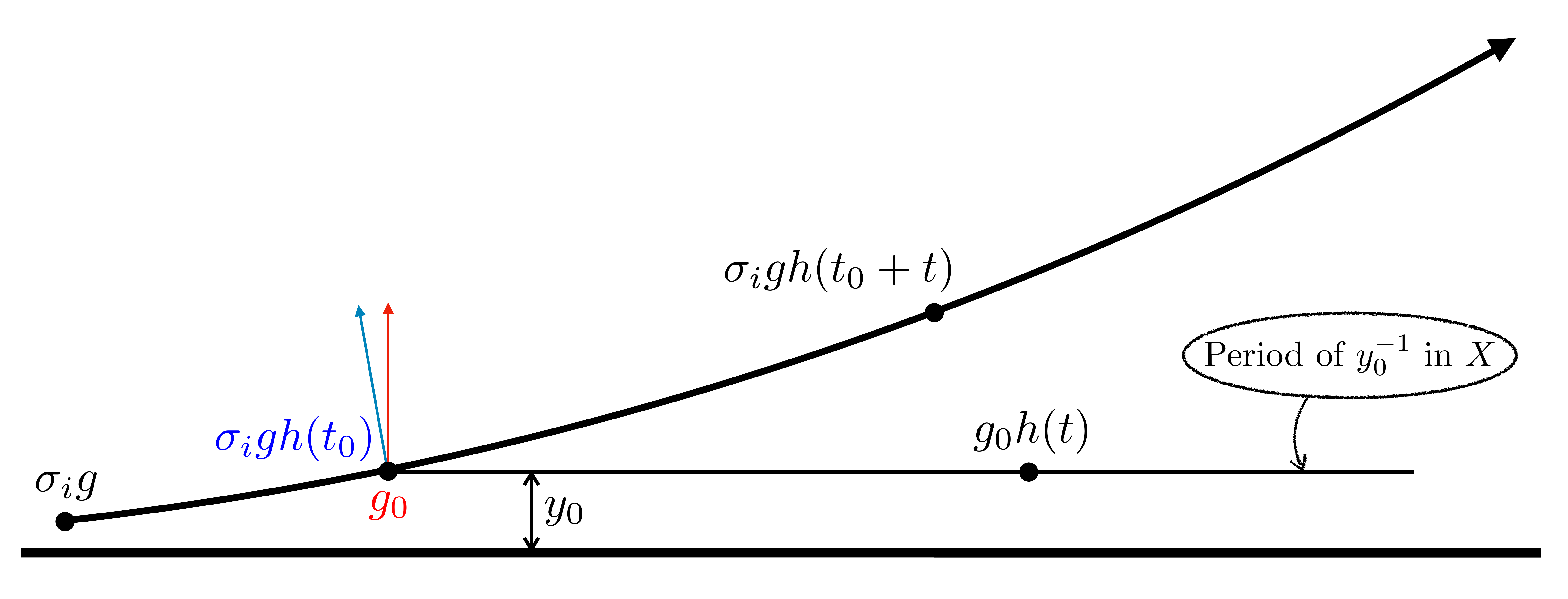} \label{figure4}\
\caption{Approximation of the segment $\{ \sigma_i gh(t_0+t), 0 \leq t \leq K\}$ with a periodic horocycle}
\end{figure}

It is thus only left to show that for all but a few exceptional $s$,  which will be the ones in an interval around $0$, this approximation is good.  To see this,  fix some $0 \leq t \leq K$ and compare 
\[
lh(s+t)=h\left(\alpha-\frac{R(s+t)}{(s+t)^2+1}\right)a\left(\frac{R}{(s+t)^2+1}\right)k(-\mathrm{arccot} \; (s+t))
\]
with 
\[
g_0h(t)=h\left(\alpha-\frac{Rs}{s^2+1}\right)a\left(\frac{R}{s^2+1}\right)h(t)=h\left(\alpha-\frac{R(s-t)}{s^2+1}\right)a\left(\frac{R}{s^2+1}\right).
\]
Firstly, note that
\[
|\mathrm{arccot}(s+t)| \ll \left| \frac{1}{s+t} \right| \leq \delta
\]
provided that $|s| \geq \delta^{-1} K$. Secondly,  note that
\[
d_G\left(a\left(\frac{R}{(s+t)^2+1}\right),  a\left(\frac{R}{s^2+1}\right) \right)=\left|\log \left(\frac{(s+t)^2+1}{s^2+1} \right) \right| \ll \delta,
\]
for any $t \leq K$ provided that $|s| \geq \delta^{-1} K$ because 
\[
\frac{d}{dt}\log \left(\frac{(s+t)^2+1}{s^2+1} \right)=\frac{2(s+t)}{(t+s)^2+1} \ll \frac{1}{|s|} \leq \delta K^{-1}.
\]
Remembering the left-invariance of the metric and using the triangle inequality,  this implies that 
\[
d_G\left(lh(s+t),h\left(\alpha-\frac{R(s+t)}{(s+t)^2+1}\right) a\left(\frac{R}{s^2+1}\right) \right) \ll \delta.
\]
Finally, 
\begin{align*}
&d_G\left(h\left(\alpha-\frac{R(s+t)}{(s+t)^2+1}\right) a\left(\frac{R}{s^2+1}\right), g_0 h(t)\right)\\
&= \frac{s^2+1}{R} d_G\left(h\left(\alpha-\frac{R(s+t)}{(s+t)^2+1}\right),  h\left(\alpha-\frac{R(s-t)}{s^2+1} \right) \right)\\
&=\frac{s^2+1}{R} \left|\frac{R(s+t)}{(s+t)^2+1}-\frac{R(s-t)}{s^2+1}\right| \\
&=\left|\frac{(s+t)(s^2+1)-(s-t)((s+t)^2+1)}{(s+t)^2+1}\right|\\
&=\left|\frac{st^2+t^3+2t}{(s+t)^2+1}\right| \ll \delta
\end{align*}
where the last inequality holds provided that $|s| \geq \delta^{-1} K^2$.  Putting everything together,  we deduce that 
\[
d_G(lh(s+t), g_0h(t)) \ll \delta
\]
provided that $|s| \geq \delta^{-1} K^2 $. We then set $\xi:=\Gamma \sigma_i^{-1} g_0$, which is a periodic horcycle with period $ y_0^{-1}$ as  
\[
\sigma_i^{-1} g_0 h(y_0^{-1})= \sigma_i^{-1} h(1)g_0 =\sigma_i^{-1} h(1) \sigma_i \sigma_i^{-1}g_0= \gamma_i \sigma_i^{-1} g_0,
\]
where we recall from the beginning of this section that $\sigma_i$ was the element such that $\sigma_i \gamma_i \sigma_i^{-1}=h(1)$ and $\gamma_i \in \Gamma$ is the unipotent element inducing the periodicity of the closed horocycles corresponding to $r_i$.
For any point $s$ such that $lh(s) \in g h([0,T])$, by definition of the fundamental period, $y_0 \geq y_T$, so that the period of the horocycle $\{\xi h(t), t \in \mathbb{R}\}$ is indeed bounded by $y_T^{-1}$.\\
Recalling $lh(s)=\sigma_i g h\left(\frac{d}{c}+s\right)$,  we can then set the exceptional interval $I_0$ to be
\[
I_0:=\left\{t \in [0, T]: \; \left|\frac{d}{c}+t \right| \leq \delta^{-1} K^2 \right\}.
\]
This assures that our estimates hold except for $t \in I_0$ and obviously,  $|I_0| \ll \delta^{-1} K^2$.\\
Regarding the second point,  we want to make sure that for any $s$ outside of an interval around $0$ we have that 
\[
\mathrm{Im(lh(s))}=\frac{R}{s^2+1} \ll \eta^{-2} y_T.
\]
Let $s_0$ be such that either $lh(s_0)= \sigma_i  g$ or $lh(s_0)=\sigma_i  gh(T)$,  depending on which of the two points minimizes the imaginary part as in the definition of $y_T$; we then have that $y_T=\frac{R}{s_0^2+1}$.  \\
Now,  the points $s$ such that $lh(s)$ lies on the horocycle orbit $\sigma_i  gh([0,T])$,  lie either in the interval $[s_0, s_0+T]$ or $[s_0-T, s_0]$,  again depending on which of the two points $\sigma_i  g, \sigma_i  gh(T)$ is minimizing.  
If $|s_0| \geq 2T$,  we have that for any such $s$ 
\[
\mathrm{Im(lh(s))} \ll \frac{R}{(|s_0|-T)^2+1}\ll y_T.
\]
If not,  we can impose $|s|>\eta T$ to assure
\[
\mathrm{Im(lh(s))} \ll \frac{R}{\eta^2T^2+1} \ll \frac{y_T T^2}{\eta^2T^2} \ll y_T \eta^{-2}.
\]
We can set $I_0$ as before, but this time with the condition $\left|s_0+t\right| \leq \eta T$.
\end{proof} 
\newpage

\section{Equidistribution of $\nu$} \label{sec::proofs}
In this section, we are going to prove Theorem \ref{main}. We will set $\theta=\frac{\beta}{40}$, where $\beta$ is the constant from Theorem \ref{venkatesh} depending only on the smallest eigenvalue of the Laplacian on $\Gamma$. In the case $\Gamma=PSL_2(\mathbb{Z})$, this makes $\frac{1}{2880}$ an admissible value for $\theta$.\\
The proof is split into different cases depending on the time parameter $T$.  To start off,  we will cover the case of good asymptotics,  where $r \gg T^{\frac{1}{20}}$, and $g_{\log T}(\xi)$ is far away from all cusps. In this case, Theorem \ref{venkatesh} is sufficient to prove good equidistribution of $\sum_{n \leq \nicefrac{T}{s}} f(\xi h(sn))$ for all $s \leq R=T^\theta$.  This immediately implies that case of  Theorem \ref{main}.
\begin{prop} \label{case1}
Let $p \in X$, $f$ such that $\Vert f \Vert=1$. Then, for all $T$ such that $r \gg T^{\frac{1}{20}}$, 
\[
\left| \frac{1}{T} \sum_{n \leq T} f(p h(n)) \nu(n)-\int f \; d\mu_X \right| \ll \frac{\log \log R}{\log R}.
\]
\end{prop}  
\begin{proof}
Assume that $\int f \; d\mu=0$, which picks up an error term of $\frac{\log \log R}{\log R}$ from the normalization of $\nu$ proven in Theorem \ref{siegelwalfisz}.  By Theorem \ref{venkatesh}, 
\[
\left| \frac{s}{T} \sum_{1 \leq sj \leq T} f(ph(sj)) \right| \leq  s^{\frac{1}{2}} r^{-\frac{\beta}{2}}.
\]
Unboxing the definition of $\nu$, we find
\begin{align*}
&\left| \frac{1}{T}\sum_{n \leq T} f(ph(n)) \nu(n) \right|= \left|\sum_{n \leq T} \sum_{\substack{e, d \leq R\\e, d |n}} \frac{\mu(e) \mu(d) f(ph(n))}{T \log R} \log\left(\frac{R}{d}\right) \log\left(\frac{R}{e}\right) \right| \\
& \leq  \sum_{e, d \leq R}\frac{\log R}{[d, e]} \left|\frac{[d, e]}{T}\sum_{n \leq \frac{T}{[e, d]}} f(ph([e, d]n)) \right| 
 \leq  \sum_{e, d \leq R}\frac{\log R}{\sqrt{[d, e]}} r^{-\frac{\beta}{2}} \\
 &\leq r^{-\frac{\beta}{2}} \sum_{m \leq R} \sum_{e, d \leq \frac{R}{m}} \frac{\log R}{\sqrt{edm}} \ll r^{-\frac{\beta}{2}} \sqrt{R}  \log R,
\end{align*}
where we ordered the terms according to their greatest common divisor $m$.
\end{proof}
In the case that $r \ll T^{\frac{1}{20}}$, we will use Lemma \ref{orbitapprox} to reduce to closed horocycles of small period and use Theorem \ref{venkatesh} together with the Siegel-Walfisz type Theorem \ref{siegelwalfisz} to conclude the proof.
\begin{proof}[Proof of Theorem \ref{main}]
Let $p \in X$ and $T$ be given such that $r \ll T^{\frac{1}{20}}$.  Assume that $\Vert f \Vert =1$ and fix some $\delta>0$ to be determined later.  Set $K:=T^{\frac{1}{3}}$.  Apply Lemma \ref{orbitapprox} to split the interval $[0,T]$ into intervals $[t_j, t_j +K]$ such that for all but on a $\delta$ proportion of them,
$\{ p h(t_j+t), t \leq K\}$ is at distance at most $\delta$ from $\{\xi_j h(t),  t \leq K\}$,  where this is a closed horocycle of period $P_j$ with $\delta^2 r \ll P_j \ll r$.  We know by Str\"ombergsson's result (\cite{strombergsson}) or from Theorem \ref{venkatesh} that
\begin{align*}
&\left| \frac{1}{T} \sum_{n \leq T} f(p h(n)) \nu(n)-\int f \; d\mu_X \right| \\
&\leq \left| \frac{1}{T} \sum_{n \leq T} f(p h(n)) \nu(n)-\frac{1}{T}\int_0^T  f(ph(t))\; dt \right|+O(r^{-\beta}) \\
&\ll 
O(r^{-\beta}+\delta)+\frac{K}{T} \sum_j \frac{1}{K} \left|  \sum_{n \leq K} f(\xi_j h(n)) \nu(n)-\int_{0}^{K} f(\xi_j h(t)) \; dt \right|.
\end{align*}
Fix some $j$ and set $y:=P_j^{-1}$.  Set $F(t):= f(\xi_j h(t y^{-1}))$, which is a $1$-periodic function and is $y^{-1}$-Lipschitz by the bounds on $f$.  It thus only remains to show that 
\[
\frac{1}{K} \sum_{n \leq K} F(yn) \nu(n)
\]
is like 
\[
\int F:=\int_0^1 F(t) dt. 
\]
We want to apply Theorem \ref{siegelwalfisz} and to do so,  we need to get from a function periodic on $[0,  1]$ to a function periodic on the integers.  To this end,  we approximate $y$ with a rational up to $R^3y^{-3}$.  That is,  we use the Dirichlet box principle to find $y^{-1} \leq q \leq R^3y^{-3}$ and $(a, q)=1$ such that 
\[
\left| y - \frac{a}{q} \right| < \frac{1}{q R^3y^{-3}}.
\]
Pick some $M$ and consider how much the function $n \mapsto F(yn)$ can diverge from a truly $q$-periodic function on an interval $\{m_0, \dots, m_0+qM\}$.  Comparing $F(y(m_0+q M))$ to $F(ym_0+aM)=F(ym_0)$ for some $m_0$,  we get that
\begin{align*}
&|F(y(m_0+q M))-F(ym_0)| \leq y^{-1} |yqM-aM| \leq \frac{My^{-1}}{R^3y^{-3}}.
\end{align*}
This is $O(y)$ provided that $q M \leq q R^3 y^{-1}$.\\
Truncate into intervals of length approximately $q R^3$; as we have just seen, the function $n \mapsto F(y n)$ is at distance $O(y)$ from a $q$-periodic function on each one.  We can thus apply Theorem \ref{siegelwalfisz} on each interval to deduce that
\[
\frac{1}{K} \sum_{n \leq K} F(yn) \nu(n)=\frac{q}{\phi(q)K} \sum_{\substack{n \leq K \\ (n,q)=1}} F(yn)+O(y)+O\left(\frac{\log \log R}{\log R}\right).
\]
To show that the sum on the right is like $\int F$,  we need one more claim.
\begin{claim} \label{smallaps} Let $\epsilon=\frac{\beta}{12}$.
\[
\left| \frac{s}{K} \sum_{sn \leq K} F(ysn) -\int F \right| \leq q^{-\epsilon}
\]
for all $s|q$ such that $s \leq q^{\epsilon}$.
\end{claim}
Before we show the claim, let us see how it allows us to conclude the proof.  
We use 
the identity $1_{m=1}=\sum_{d|m} \mu(d)$ to find
\begin{align*}
\sum_{\substack{n \leq K \\ (n, q)=1}} F(yn) &=\sum_{n \leq K} \sum_{d|(n, q)} \mu(d) F(yn)\\
&=\sum_{d|q} \mu(d) \sum_{\substack{n \leq K \\ d|n}} F(yn)=\sum_{d|q} \mu(d) \sum_{nd \leq K}F(ydn).
\end{align*}
Decomposing the sum and using Claim \ref{smallaps}, we see that
\begin{align*}
\frac{q}{\phi(q)K} \sum_{\substack{n \leq K \\ (n, q)=1}} F(yn) &\leq \frac{q}{\phi(q)} \sum_{\substack{d|q\\d<q^\epsilon}} \frac{1}{d} \left| \frac{d}{K} \sum_{dn \leq K}F(ydn)\right|+\frac{q}{\phi(q)} \sum_{\substack{d|q\\d \geq q^\epsilon}} \frac{\Vert F \Vert_\infty}{d}\\
& \leq 2 \frac{q^{1-\epsilon} \tau(q)}{\phi(q)}.
\end{align*}
By standard asymptotics of $\phi$ and $\tau$ (see for example \cite{iwaniec}), the right-hand side is 
 \[
 O(q^{-\frac{6\epsilon}{7}})=O(y^{\frac{\beta}{14}})=O\left(\delta^{-\frac{\beta}{7}} r^{\frac{\beta}{14}} \right). 
 \]
We can choose $\delta=r^{-\frac{1}{5}}$ to get the desired conclusion.  It thus only remains to show Claim \ref{smallaps}.\\
To prove Claim \ref{smallaps}, we divide into two cases. Firstly,  in the case that $q \leq y^{-3}$,  we apply Str\"ombergsson's result (or Theorem \ref{venkatesh}) to the periodic horocycle $\xi_j h(t)$, for which  $r(\xi_j, T)=y^{-1}$ for any $T$, to see
\[
\int F=y \int_{0 \leq t \leq y^{-1}} f(\xi_j h(t)) dt=\int\;  f \; d\mu_X+O(y^\beta).
\]
Using this and applying Theorem \ref{venkatesh} to the same periodic horocycle piece, we see
\[
\forall s \leq y^{\frac{\beta}{4}}: \quad \left| \frac{s}{K} \sum_{1 \leq j \leq \nicefrac{K}{s}} F(ysj)-\int F  \right| \ll s^{\frac{1}{2}} y^{\frac{\beta}{2}} \ll y^{\frac{\beta}{4}}.
\]
As $q^{\epsilon} \leq y^{\frac{\beta}{4}}$,  we deduce Claim \ref{smallaps} in this case.\\
For the second case,  assume that $q \geq y^{-3}$.  Roughly speaking,  in this case there are so many distinct points of the form $sn\frac{a}{q}$ in the interval $[0, 1]$ that they cannot help being dense enough to approximate $\int F$ by force.  \\
We split $[0, K]$ into intervals $I$ of length $q R^3 y^{-1}$.  Fix $s \leq q^{\epsilon}$ and set $q^\prime:=\nicefrac{q}{s}$.  Fix an interval $I$ and call its left endpoint $s t_0$.  We note that as for any $n$ such that $sn \in I$,
\begin{align*}
\left|F(ysn)-F\left(yst_0+\frac{sa}{q}(n-t_0)\right)\right| &\leq y^{-1} \left|y-\frac{a}{q}\right| |sn-st_0|\\
& \leq \frac{y^{-1}|I|}{q R^3y^{-3}}\leq y;
\end{align*}
set $x_0:=st_0\left(y-\frac{a}{q}\right)$ and note that then
\begin{align*}
&\frac{s}{|I|} \sum_{sn \in I} F(ysn) = O(y)+\frac{s}{|I|} \sum_{sn \in I} F\left(x_0+n\frac{as}{q}\right)\\
&=O(y)+ O\left(\frac{s q^\prime}{|I|} \right)+\frac{s}{q} \sum_{n \leq q^\prime} F\left(x_0+n\frac{a}{q^\prime}\right),
\end{align*}
where we use in the second line that the function $F(x_0+\frac{san}{q})$  is $q^\prime$ periodic in $n$.  The number $a$ is coprime to $q^\prime$,  so it plays no role and can be dropped; we then only have to evaluate
\[
\frac{1}{q^\prime} \sum_{n \leq q^\prime} F\left(x_0+\frac{n}{q^\prime}\right).
\]
But for any $t \in (0,1)$ and any $n$,
\[
\left| F\left(x_0+\frac{n}{q^\prime}\right)-F\left(x_0+\frac{n+t}{q^\prime}\right) \right| \leq y^{-1} \frac{1}{q^\prime} \leq y^{-1} y^{3(1-\epsilon)} \leq y,
\]
which implies that
\begin{align*}
&\frac{1}{q^\prime} \sum_{n \leq q^\prime} F\left(x_0+\frac{n}{q^\prime}\right)=O(y)+\frac{1}{q^\prime} \sum_{n \leq q^\prime} \int_0^1 F\left(x_0+\frac{n+t}{q^\prime}\right) \; dt\\
&=O(y)+\frac{1}{q^{\prime}} \int_0^{q^{\prime}} F\left(\frac{t}{q^\prime} \right) \; dt=O(y)+\int F
\end{align*}
This shows Claim \ref{smallaps} also in the second case, which, as we have seen, concludes the proof of Theorem \ref{main}.
\end{proof}
\newpage

\printbibliography

\end{document}